# Cascadic Tensor Multigrid Method and Economic Cascadic Tensor Multigrid Method for Image Restoration Problems

**Ziqi Yan, Chenliang Li*, Yuhan Chen**

School of Mathematics and Computing Science, Center for Applied Mathematics of Guangxi (GUET),
Guilin University of Electronics Technology, Guilin 541004, China
*Corresponding author: chenli@guet.edu.cn

**Abstract** A cascadic tensor multigrid method and an economic cascadic tensor multigrid method is presented for solving the image restoration models. The methods use quadratic interpolation as prolongation operator to provide more accurate initial values for the next fine grid level, and constructs a preserving-edge-denoising operator to obtain better edges and remove noise. The experimental results show that the new methods not only improves computational efficiency but also achieve better restoration quality.

*Keywords*: Image restoration, Cascadic tensor multigrid methods, Prolongation operator, Preserving-edge-denoising operator

## 1. Introduction

As a common information carrier, in the process of formation, transmission, recording, and processing, an image is easy to be affected by the equipment, environment, and other impacts of blurring, distortion, and so on. It is very important to recover effectively the blurred image. In recent years, image retoration problem has received extensive attention from scholars, and many models and algorithms have been proposed, such as super-resolution reconstruction[1], image denoising [2], image restoration [3], and so on.

In addition to the traditional vector-based image restoration model, some scholars use tensor methods to solve the image restoration problem. [4] constructs a tensor representation-based degradation model that applies tensor CP decomposition to color image and video recovery. [5] proposes a new tensor degeneracy model for color image and video recovery. The tensor degeneracy model can be recovered by solving the tensor equation under the Einstein-product.

Consider the following tensor color image model [6],

$$G(x, y, z) = F(x, y, z) * S(x, y, z) + N(x, y, z) \quad (1)$$

where $G(x, y, z)$ denotes the blurred image, $F(x, y, z)$ denotes the original image, $*$ denotes the 3D convolution operator, $S(x, y, z)$ denotes the 3D Gaussian function, and $N(x, y, z)$ denotes the noise function. The corresponding three-dimensional model of image can be ritten as the following tensor equation,

$$\mathcal{G} = \mathcal{T} *_3 \mathcal{F} + \mathcal{N}, \quad (2)$$

where $\mathcal{T} \in \mathbb{R}^{I_1 \times I_2 \times 3 \times I_1 \times I_2 \times 3}$ denotes a tensor blurred operator, $\mathcal{G} \in \mathbb{R}^{I_1 \times I_2 \times 3}$ denotes tensor blurred images, $\mathcal{F} \in \mathbb{R}^{I_1 \times I_2 \times 3}$ denotes the image to be recovered, $\mathcal{N}$ denotes random noise, $*$ denotes the Einstein-product.

The tensor Einstein-product [7] is defined as follows: let $\mathcal{T} \in \mathbb{R}^{I_1 \times I_2 \times \ldots \times I_L \times K_1 \times K_2 \times \ldots \times K_N}$, $\mathcal{F} \in \mathbb{R}^{K_1 \times K_2 \times \ldots \times K_N \times J_1 \times J_2 \times \ldots \times J_M}$. Its elements are defined as

$$(\mathcal{T} *_N \mathcal{F})_{i_1 \ldots i_L j_1 \ldots j_M} = \sum_{k_1, \ldots, k_N} t_{i_1 \ldots i_L k_1 \ldots k_N} f_{k_1 \ldots k_N j_1 \ldots j_M} \quad (3)$$

In order to solve the tensor equation (3) under the Einstein-product, scholars have proposed a number of effective methods. Brazell et al. [8] established a high order dual Conjugate Gradient method (CG) and a high order Jacobi methods. Wang et al. [9] proposed CG for solving the tensor equation under the Einstein-product and proved the finite termination of this method. Huang et al. [10] proposed Conjugate Residual method based on a tensor format (CR-BTF), Conjugate Gradient-Schmidt method based on a tensor format (CGS-BTF), and Bi-Conjugate Gradient method based on a tensor format (BiCG-BTF), etc. The convergence analysis of the methods is given, and the numerical results verify the feasibility of the methods.

Solving the tensor equation under the Einstein-product has some limitations when dealing with high-resolution color images, which require a large amount of memory. The cascadic resolution method [10] is a good image recovery algorithm.

Cascadic multigrid method is a simple and efficient method, please see [12] and therein references. Based on the super-convergence and extrapolation methods, [13] constructs a better interpolation operator, to get a good initial value for the next iteration to speed up the convergence of the algorithm. [11] proposes a class of nonlinear preserving-edge-denoising operator as the prolongation operator, which can reduce the ringing effect and recover good image edges. [14] uses a new quadratic interpolation method that can provide better initial values for fine grid levels. [15] further investigates the selection of blurred operator and smoothing methods and shows that the cascadic multigrid method has good efficiency for image recovery. An economic cascadic multigrid method is proposed in [16]. Compared with the traditional cascadic multigrid method, this method proposes a new formula for the number of iterations, which can obtain the same numerical solution with fewer iterations.

In this paper, we construct the tensor form of the edge-preserving-denoising-operator and the prolongation operator and propose cascadic tensor multigrid method and an economic cascadic tensor multigrid algorithm for the tensor color image model (1). Finally, some experiments show that our algorithm is more efficient.

The structure of this paper is as follows. Section 2, introduces the basic definition of tensor Einstein-product and introduces the basic framework of the economic cascadic tensor multigrid method. Section 3, introduces the prolongation operator. Section 4, introduces edge-preserving-denoising-operator. In Section 5, some numerical results are given to illustrate the effectiveness of the new method. In Section 6, a summary is presented.

## 2. Preliminaries

The cascadic multigrid method is an efficient one-way multigrid method [17]. Compared with linear interpolation, quadratic interpolation can provide better initial values for the fine grid level. [13] presents the following cascadic multigrid method.

**Algorithm 1: Cascadic Multigrid Method (CMG)** [13]
**Step 1**: When $j=1$, exact solution on the coarsest grid level $\mathbf{u}_1^1 = \mathbf{u}_1^* = A_1^{-1} f_1$;
**Step 2**: When $j = 2, 3, \cdots, L$, iterative calculation:

(1) Prolongation: $\mathbf{u}_j^1 = I_j^2 \mathbf{u}_{j-1}^*$;

(2) Smoothing: $\mathbf{u}_j^{m_j} = G_j^{m_j}\left(\mathbf{A}_j, \mathbf{f}_j, \mathbf{u}_j^1\right)$;

(3) $\mathbf{u}_j^* = \mathbf{u}_j^{m_j}$.

Where $I_j^2(\cdot)$ denotes the quadratic interpolation operator; $G_j^{m_j}(\cdot)$ denotes the smoothing of $m_j$ times in the $j$ th level, which satisfies $m_j = [m_* l^2]$, where $m_* \geq 1$, and can be chosen as a smoothing operator by a method such as conjugate gradient method.

Xu, Shi and Huang [16] present an economical cascadic multigrid method, which requires less work operations on the each level, through control of the iteration numbers on the each level to preserve the accuracy without over iterations.

**Algorithm 2: Economic Cascadic Multigrid Method (ECMG)** [16]
**Step 1**: When $j=1$, exact solution on the coarsest grid level $\mathbf{u}_1^1 = \mathbf{u}_1^* = A_1^{-1} f_1$;
**Step 2**: When $j = 2, 3, \cdots, L$, iterative calculation:

(1) Prolongation: $\mathbf{u}_j^1 = I_j^2 \mathbf{u}_{j-1}^*$;

(2) Smoothing: $\mathbf{u}_j^{m_j} = G_j^{m_j}\left(\mathbf{A}_j, \mathbf{f}_j, \mathbf{u}_j^1\right)$;

(3) $\mathbf{u}_j^* = \mathbf{u}_j^{m_j}$.

The number of iterations $m_j$ in the algorithm is chosen as follows.

(1) If $l > L_0$, then
$$m_j = \left[ m_0 (L - L_0)^2 \beta^{L-l} \right],$$

(2) If $l \leq L_0$, then
$$m_j = \left[ m_*^{\frac{1}{2}} \left(L - (2 - \varepsilon_0)l\right) h_l^{-2} \right],$$

Where $l = 1, 2, \ldots, L$, $L_0 = L/2$, $m_0 = 1$, $\beta = 4$, $m_* = 1$, $\varepsilon_0 = 1/2$, $h_l = (1/2)^{l-1}$. Here, $[t]$ denotes the smallest positive integer not less than $t$.

In [18], a cascadic multigrid method for solving the image recovery problem is proposed. When dealing with color images, it is necessary to extend the method to three dimensions. The use of tensors is more advantageous for the processing of multidimensional data.

In [19], Chen and Li presented the tensor multigrid method to the Sylvester tensor equations, maintaining accuracy while potentially reducing time consumption. So we combine the ideas of tensor multigrid and cascadic multigrid methods, propose a cascadic tensor multigrid algorithm to solve the tensor equations under the Einstein-product to recover images.

**Algorithm 3: Cascadic Tensor Multigrid Method (CTMG)**
**Step 1**: When $l = 1$, the coarsest grid level is solved exactly for $\mathcal{F}_1^\delta$: $\mathcal{T}_1 *_3 \mathcal{F}_1^\delta = \mathcal{G}_1$,
**Step 2**: When $l = 2, 3, \ldots, L$, Iterative calculation:

(1) Prolongation: $\mathcal{F}_{l,k}^0 := P_l(\mathcal{F}_{l-1,k}^\delta)$;

(2) Preserving-edge-denoising: $\mathcal{F}_{l,k}^D := D_l(\mathcal{F}_{l,k}^0)$;

(3) Smoothing: $\mathcal{F}_l^\delta = G_j^{m_j}\left(\mathcal{T}_l, \mathcal{G}_l, \mathcal{F}_l^\delta\right)$;

**Step 3**: $\mathcal{F}^\delta = \mathcal{F}_L^\delta$.

where $L$ denotes the maximum number of level of the grid; $P_l(\cdot)$ denotes prolongation operator; $D_l(\cdot)$ denotes preserving-edge-denoising operator. The construction of $P_l(\cdot)$ and $D_l(\cdot)$ are introduced in Section 3 and Section 4. $G_j^{m_j}(\cdot)$ denotes cascadic multigrid method, where $m_j$ denotes smoothing $m_j$ times at the $j$ th level and satisfies $m_j = [m_* l^2]$, where $m_* \geq 1$. Conjugate Residual method based on tensor form (CR-BTF), Conjugate Gradient-Schmidt method (CGS-BTF), or Bi-Conjugate Gradient method (BiCG-BTF) etc. [10] can be used as the smoothing operator.

Similarly we can obtain the following economic cascadic tensor multigrid method.

**Algorithm 4: Economic cascadic tensor multigrid approach (ECTMG)**

**Step 1**: When $l = 1$, the coarsest grid level is solved exactly for $\mathcal{F}_1^\delta$: $\mathcal{T}_1 *_3 \mathcal{F}_1^\delta = \mathcal{G}_1$,

**Step 2**: When $l = 2, 3, \ldots, L$, iterative computation:

(1) Prolongation: $\mathcal{F}_{l,k}^0 := P_l \cdot \mathcal{F}_{l-1,k}^\delta$;

(2) Preserving-edge-denoising: $\mathcal{F}_{l,k}^D := D_l \cdot \mathcal{F}_{l,k}^0$;

(3) Smoothing: $\mathcal{F}_l^\delta = G_j^{m_j}\left(\mathcal{T}_l, \mathcal{G}_l, \mathcal{F}_l^\delta\right)$,.

**Step 3**: $\mathcal{F}^\delta = \mathcal{F}_L^\delta$.

The number of iterations $m_j$ in the algorithm is chosen as follows.

(1) If $l > L_0$, then
$$m_j = \left[m_0(L - L_0)^2 \beta^{L-l}\right],$$

(2) If $l \leq L_0$, then
$$m_j = \left[m_*^{\frac{1}{2}}(L - (2 - \varepsilon_0)l) h_l^{-2}\right],$$

where, $\beta = 4$ $L_0 = L/2$, $m_0 = 1$, $\beta = 4$, $m_* = 1$, $\varepsilon_0 = 1/2$, $h_l = (1/2)^{l-1}$. Here, $[t]$ denotes the smallest positive integer not less than $t$. CR-BTF, CGS-BTF, BiCG-BTF [10] can be used as the smoothing operator.

## 3. Prolongation Operator

To improve the convergence speed of the algorithm, it is necessary to construct an effective prolongation operator to obtain a better initial value for the next level.

The image can be regarded as a tensor composed of R, G and B channels. In each channel, the image will be prolonged by the quadratic interpolation operator. The construction of the specific quadratic interpolation operator is introduced in [18].

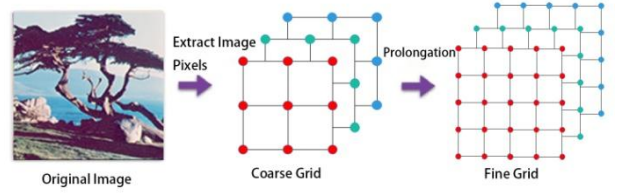

**Figure 1.** Processing of three channels of an image

In [18], as shown in (Figure 1), the hollow circles represent pixels at the $l_2$ grid level, the solid black dots represent pixel values at the $l_1$ grid level, and the rectangular boxes represent pixel values at the $l$ grid level. The image values at the other color nodes are computed by quadratic interpolation.

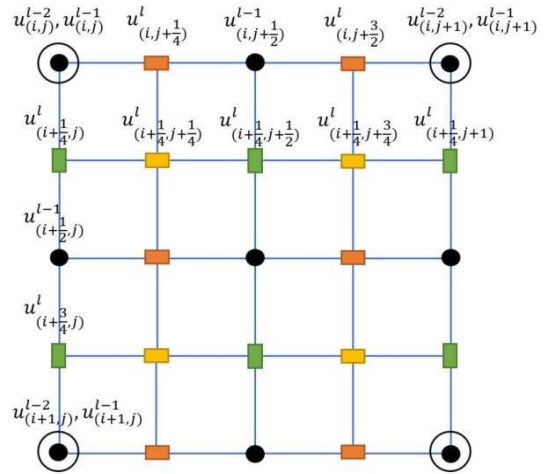

**Figure 2.** Quadratic interpolation methods for the two-dimensional case [18]

Let $f_{(i,j,k)}^l$ denote the pixel values in $(i, j, k)$ on the $l$ level, $k = 1, 2, 3$ denotes R,G,B channels respectively.

$$f_{(i,j+1/4,k)}^l = \frac{1}{16}\left[\left(9f_{(i,j,k)}^{l-1} + 12f_{(i,j+1/2,k)}^{l-1} - f_{(i,j+1,k)}^{l-1}\right) - \left(3f_{(i,j,k)}^{l-2} + f_{(i,j+1,k)}^{l-2}\right)\right],$$

$$f_{(i,j+3/4,k)}^l = \frac{1}{16}\left[\left(9f_{(i,j+1,k)}^{l-1} + 12f_{(i,j+1/2,k)}^{l-1} - f_{(i,j,k)}^{l-1}\right) - \left(3f_{(i,j+1,k)}^{l-2} + f_{(i,j,k)}^{l-2}\right)\right],$$

$$f_{(i+1/4,j,k)}^l = \frac{1}{16}\left[\left(9f_{(i,j,k)}^{l-1} + 12f_{(i+1/2,j,k)}^{l-1} - f_{(i+1,j,k)}^{l-1}\right) - \left(3f_{(i,j,k)}^{l-2} + f_{(i+1,j,k)}^{l-2}\right)\right],$$

$$f_{(i+3/4,j,k)}^l = \frac{1}{16}\left[\left(9f_{(i+1,j,k)}^{l-1} + 12f_{(i+1/2,j,k)}^{l-1} - f_{(i,j,k)}^{l-1}\right) - \left(3f_{(i+1,j,k)}^{l-2} + f_{(i,j,k)}^{l-2}\right)\right],$$

$$f_{(i+1/4,j+1/4,k)}^l = \frac{3}{8}f_{(i+1/4,j,k)}^l + \frac{3}{4}f_{(i+1/4,j+1/2,k)}^l - \frac{1}{8}f_{(i+1/4,j+1,k)}^l,$$

$$f_{(i+1/4,j+3/4,k)}^l = -\frac{1}{8}f_{(i+1/4,j,k)}^l + \frac{3}{4}f_{(i+1/4,j+1/2,k)}^l + \frac{3}{8}f_{(i+1/4,j+1,k)}^l.$$

## 4. Tensor Preserving-Edge-Denoising Operator

The color image is a composite image consisting of three "gray matter" maps representing the R, G, and B channels. [18] describes the preserving-edge-denoising operator in matrix form. The preserving-edge-denoising operator is performed for images under R, G, and B channels respectively. In each channel, the basic principle is the $P-M$ equation [16],

$$\begin{cases} \dfrac{\partial u}{\partial t} = \mathrm{div}[g(|\nabla u|)\nabla u] \\ u(x,y,t)|_{t=0} = u_0(x,y) \end{cases}, (x,y)\in R^2, \quad (4)$$

where $u(x,y,t)$ denotes the evolving image at time $t$, and $g(|\nabla u|)\in[0,1]$ is the edge stopping function or a function of the diffusion coefficient that serves to make the image smooth inside the edges and therefore tends to 1. The diffusion coefficients are as follows.

$$g(|\nabla u|) = \dfrac{1}{1+\left(\dfrac{|\nabla u|}{k}\right)^2}, \quad (5)$$

where $k$ denotes the threshold value, which can be pre-set or changed with the result of each iteration of the image, and it is related to the variance of the noise. Let $c(x,y)=g(|\nabla u|)$, the image preserving-edge-denoising process can be expressed in matrix form as

$$\mathbf{u}^{n+1} = \mathbf{u}^n + \tau \mathbf{A}\left(\mathbf{u}^n\right)\mathbf{u}^n, \quad (6)$$

Where $\mathbf{u}^n$ and $\mathbf{u}^{n+1}$ denote the $n$ and $n+1$ moment images, respectively. $\mathbf{A}\left(\mathbf{u}^n\right)$ denotes the $N^2\times N^2$ matrix.

We perform (6) 10 times to achieve the effect of preserving image edges and removing noise.

## 5. Numerical Results

The image restore model can be expressed as the following tensor, i.e.

$$\mathcal{T}*_3 \mathcal{F} = \mathcal{G} + \mathcal{N}, \quad (1)$$

The blurred operator tensor $\mathcal{T}_l \in \mathbb{R}^{n_1\times n_2\times 3\times n_1\times n_2\times 3}$ is a sixth-order Toeplitz tensor [20-22] obtained from $\mathcal{S}$, that is,

$$\begin{aligned}&\mathcal{T}(i_1,i_2,i_3,j_1,j_2,j_3)\\&=\mathcal{S}(j_1-i_1,j_2-i_2,j_3-i_3),\quad 1\leq i_k,j_k\leq n, 1\leq k\leq m\end{aligned} \quad (2)$$

where $S$ denotes a three-dimensional Gaussian function.

$$S(i,j,k) = \dfrac{1}{(\sqrt{2\pi}\sigma)^3}\exp\left(-\dfrac{i^2+j^2+k^2}{2\sigma^2}\right) \quad (9)$$

In our numerical experiments, let $\sigma=0.7,0.8,0.9$, we obtain different blurred operators and get different blurred images. And we add random noise $\mathcal{N}=0.001*rand(100,100,3)$ to get the blurred and noise image $\mathcal{G}$. The correspondly image are shown in (Figure 3).

Numerical evaluation indexes of image restoration are relative error (denoted by "RE") and the peak signal to noise ratio (denoted by "PSNR"), which are widely used in visual data restoration tasks.

$$\mathrm{RE} = \dfrac{\|\mathcal{G}-\mathcal{F}\|_F}{\|\mathcal{F}\|_F}, \quad \mathrm{PSNR} = 10\log_{10}\left(\dfrac{n_1n_2n_3\mathcal{F}_{\max}^2}{\|\mathcal{G}-\mathcal{F}\|_F^2}\right),$$

where $\mathcal{F}$ denotes the original image, $\mathcal{G}$ denotes the recovered image, and $n_1,n_2,n_3$ are the dimensions of the tensor $\mathcal{F}$. $\mathcal{F}_{\max}$ denotes the maximum pixel value of the tensor $\mathcal{F}$ of the original image. The size of the image is $128\times 128$.

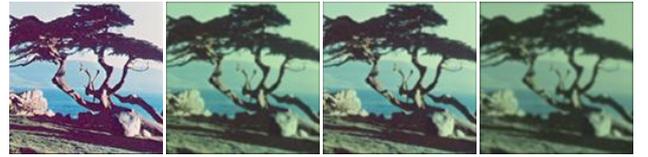

a.Original image    b. $\sigma=0.7$    c. $\sigma=0.8$    d. $\sigma=0.9$

**Figure 3.** Original image and blurred image

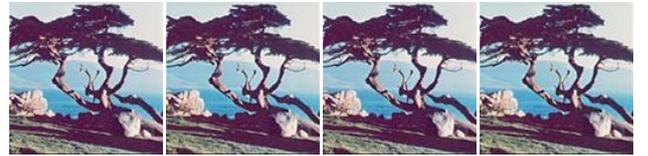

a.Original image    b.BICG-BTF PSNR=27.557    c.CGS-BTF PSNR=27.372    d.CR-BTF PSNR=27.425

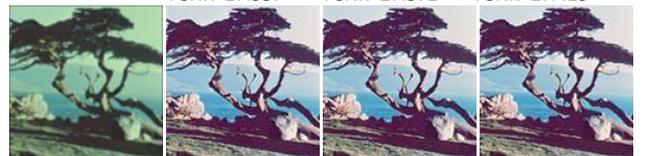

e.Blurred image    f.CTMG-BICG PSNR=29.178    g.CTMG-CGS PSNR=29.345    h.CTMG-CR PSNR=29.634

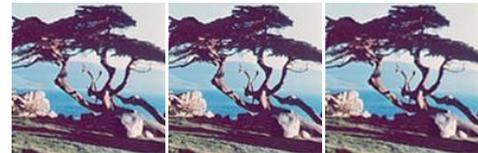

i. ECTMG-BICG PSNR=28.950    j.ECTMG-CGS PSNR=29.524    k.ECTMG-CR PSNR=29.810

**Figure 4.** Recovered Image of tree

(Figure 4) and (Figure 5) show that CTMG methods, ECTMG methods are better than BICG-BTF, CGS-BTF and CR-BTF.

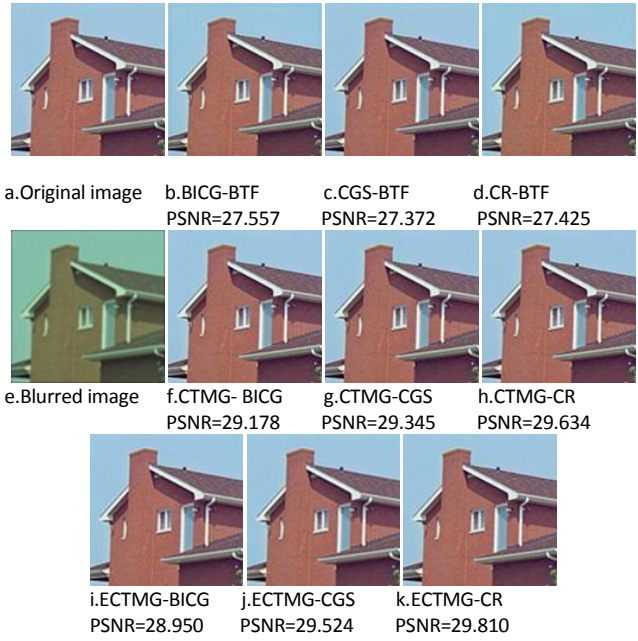

a. Original image  b. BICG-BTF  c. CGS-BTF  d. CR-BTF
PSNR=27.557  PSNR=27.372  PSNR=27.425

e. Blurred image  f. CTMG-BICG  g. CTMG-CGS  h. CTMG-CR
PSNR=29.178  PSNR=29.345  PSNR=29.634

i. ECTMG-BICG  j. ECTMG-CGS  k. ECTMG-CR
PSNR=28.950  PSNR=29.524  PSNR=29.810

**Figure 5.** Recovered Image of house

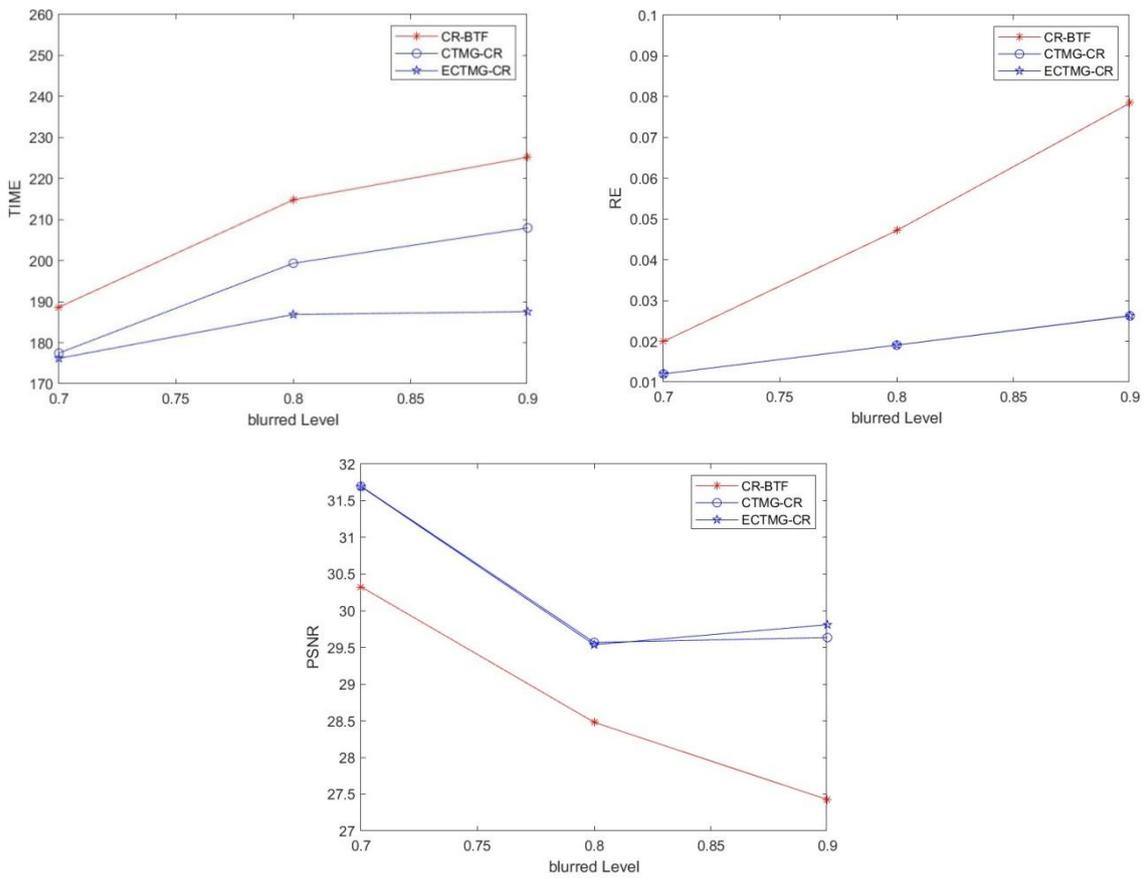

**Figure 6.** Comparison of BICG-BTF, CTMG-BICG and ECTMG-BICG for tree(128 × 128).

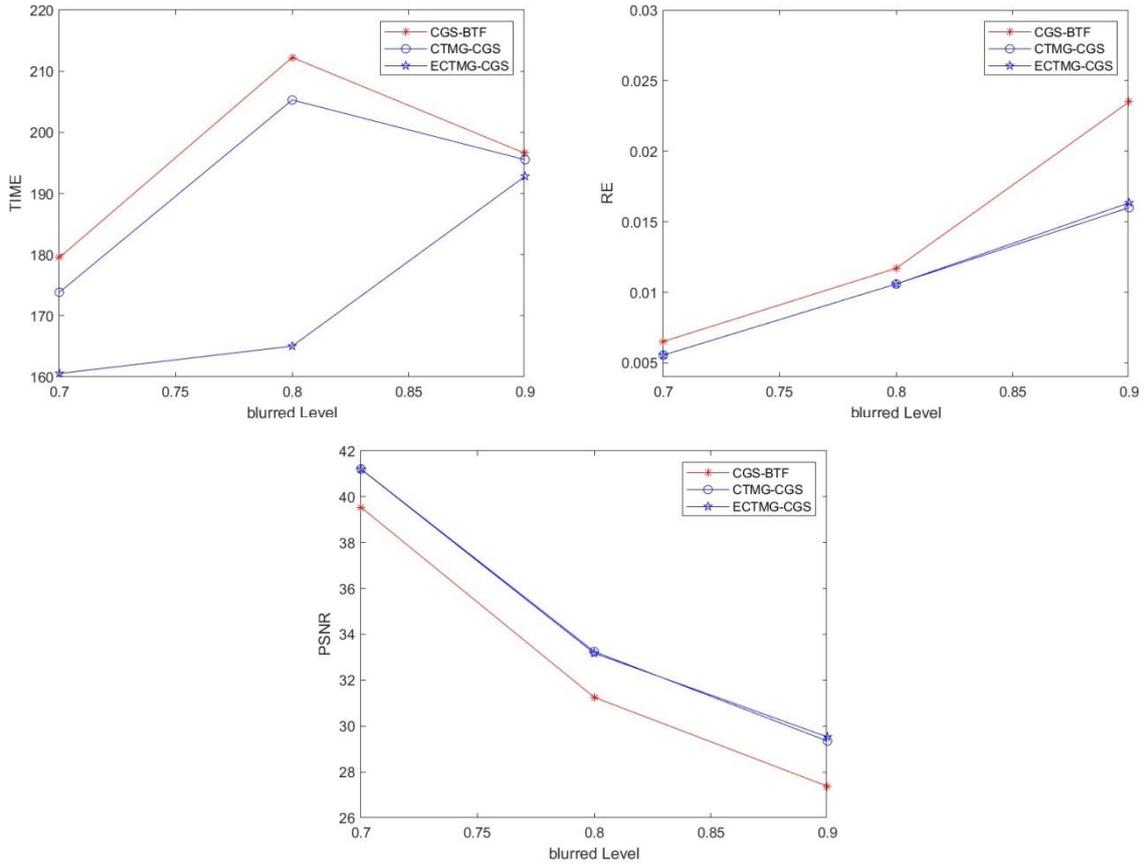

**Figure 7.** Comparison of CGS-BTF, CTMG-CGS and ECTMG-CGS for tree(128 × 128).

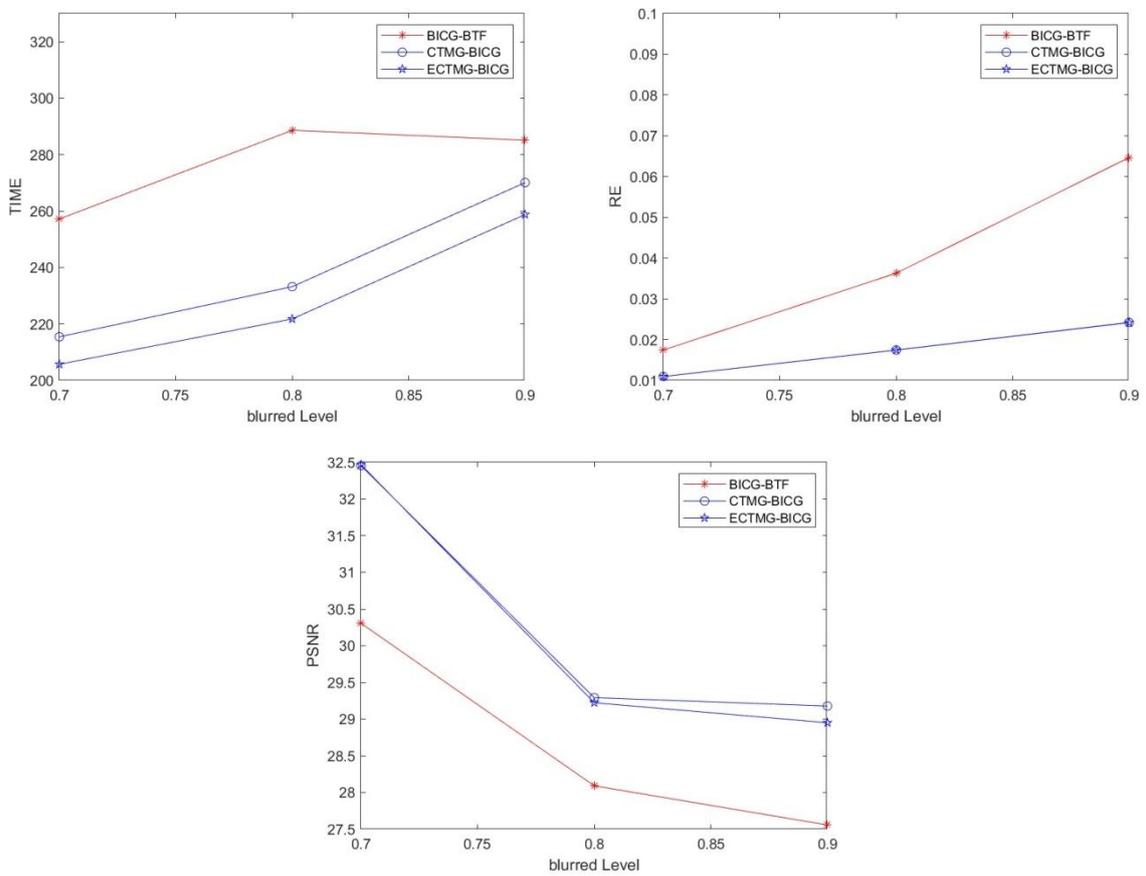

**Figure 8.** Comparison of CR-BTF, CTMG-CR and ECTMG-CR for tree(128 × 128).

**Table 1. Numerical results of BICG-BTF, CTMG-BICG and ECTMG-BICG for tree(128×128)**

| Algorithm | BICG-BTF | | | CTMG-BICG | | | ECTMG-BICG | | |
|---|---|---|---|---|---|---|---|---|---|
| Noise Fuzzy Level | CPU(s) | PSNR | RE | CPU(s) | PSNR | RE | CPU(s) | PSNR | RE |
| $\sigma = 0.7$ | 257.180 | 30.305 | 0.017 | 215.419 | 32.454 | 0.011 | 205.697 | 32.467 | 0.011 |
| $\sigma = 0.8$ | 288.609 | 28.090 | 0.036 | 233.205 | 29.292 | 0.017 | 221.760 | 29.223 | 0.017 |
| $\sigma = 0.9$ | 285.113 | 27.557 | 0.065 | 270.083 | 29.178 | 0.024 | 258.762 | 28.950 | 0.024 |

**Table 2. Numerical results of CGS-BTF, CTMG-CGS and ECTMG-CGS for tree(128×128)**

| Algorithm | CGS-BTF | | | CTMG-CGS | | | ECTMG-CGS | | |
|---|---|---|---|---|---|---|---|---|---|
| Noise Fuzzy Level | CPU(s) | PSNR | RE | CPU(s) | PSNR | RE | CPU(s) | PSNR | RE |
| $\sigma = 0.7$ | 179.528 | 39.542 | 0.007 | 173.844 | 41.209 | 0.006 | 160.567 | 41.201 | 0.006 |
| $\sigma = 0.8$ | 212.255 | 31.252 | 0.012 | 205.305 | 33.247 | 0.011 | 165.038 | 33.176 | 0.011 |
| $\sigma = 0.9$ | 196.570 | 27.372 | 0.023 | 195.526 | 29.345 | 0.016 | 192.809 | 29.524 | 0.016 |

**Table 3. Numerical results of CGS-BTF, CTMG-CGS and ECTMG-CGS for tree(128×128)**

| Algorithm | CR-BTF | | | CTMG-CR | | | ECTMG-CR | | |
|---|---|---|---|---|---|---|---|---|---|
| Noise Fuzzy Level | CPU(s) | PSNR | RE | CPU(s) | PSNR | RE | CPU(s) | PSNR | RE |
| $\sigma = 0.7$ | 188.647 | 30.326 | 0.020 | 177.469 | 31.696 | 0.012 | 176.186 | 31.692 | 0.012 |
| $\sigma = 0.8$ | 214.823 | 28.479 | 0.047 | 199.358 | 29.569 | 0.019 | 186.866 | 29.537 | 0.019 |
| $\sigma = 0.9$ | 225.217 | 27.425 | 0.078 | 207.971 | 29.634 | 0.026 | 187.576 | 29.810 | 0.026 |

**Table 4. Numerical results of BICG-BTF, CTMG-BICG and ECTMG-BICG for house(128×128)**

| Algorithm | BICG-BTF | | | CTMG-BICG | | | ECTMG-BICG | | |
|---|---|---|---|---|---|---|---|---|---|
| Noise Fuzzy Level | CPU(s) | PSNR | RE | CPU(s) | PSNR | RE | CPU(s) | PSNR | RE |
| $\sigma = 0.7$ | 246.649 | 32.461 | 0.026 | 228.445 | 34.564 | 0.019 | 233.878 | 34.546 | 0.019 |
| $\sigma = 0.8$ | 278.097 | 28.859 | 0.047 | 261.187 | 31.472 | 0.032 | 223.890 | 31.503 | 0.032 |
| $\sigma = 0.9$ | 299.900 | 26.707 | 0.045 | 287.142 | 28.100 | 0.035 | 275.066 | 28.806 | 0.035 |

**Table 5. Numerical results of CGS-BTF, CTMG-CGS and ECTMG-CGS for house(128×128)**

| Algorithm | CGS-BTF | | | CTMG-CGS | | | ECTMG-CGS | | |
|---|---|---|---|---|---|---|---|---|---|
| Noise Fuzzy Level | CPU(s) | PSNR | RE | CPU(s) | PSNR | RE | CPU(s) | PSNR | RE |
| $\sigma = 0.7$ | 177.493 | 38.405 | 0.011 | 170.162 | 41.156 | 0.009 | 158.921 | 41.172 | 0.009 |
| $\sigma = 0.8$ | 176.322 | 33.982 | 0.021 | 173.163 | 34.767 | 0.018 | 164.681 | 34.744 | 0.018 |
| $\sigma = 0.9$ | 201.217 | 29.775 | 0.040 | 168.076 | 31.887 | 0.030 | 162.674 | 31.557 | 0.030 |

**Table 6. Numerical results of CR-BTF, CTMG-CR and ECTMG-CR for house(128×128)**

| Algorithm | CR-BTF | | | CTMG-CR | | | ECTMG-CR | | |
|---|---|---|---|---|---|---|---|---|---|
| Noise Fuzzy Level | CPU(s) | PSNR | RE | CPU(s) | PSNR | RE | CPU(s) | PSNR | RE |
| $\sigma = 0.7$ | 237.442 | 32.106 | 0.029 | 166.963 | 33.906 | 0.022 | 190.228 | 33.906 | 0.022 |
| $\sigma = 0.8$ | 218.509 | 28.738 | 0.057 | 206.352 | 30.457 | 0.036 | 167.286 | 30.537 | 0.036 |
| $\sigma = 0.9$ | 256.396 | 26.509 | 0.050 | 238.227 | 27.658 | 0.034 | 212.165 | 27.918 | 0.034 |

From (Table 1) to (Table 6) and (Figure 6) to (Figure 8), the economic cascadic tensor multigrid, approach achieve higher PSNR and cost less CPUs. ECTMG-CGS performs best with highest PSNR and least CPUs among these methods.

## 6. Summaries

In this paper, by combining tensor and cascadic multigrid, we present a tensor cascadic multigrid method and an economic tensor cascadic tensor multigrid method. The numerical results show that, compared with the existing methods, methods recover images with higher PSNR and smaller relative error, while reducing the CPU time. For the image restoration problems, CTMG and ECTMG perform better whether PSNR, RE or CPUs.

## ACKNOWLEDGMENTS

This work was supported by Natural Science Foundation of China (12161027), Guangxi Natural Science Foundation (2020GXNSFAA159143).